\documentclass[11pt]{article}
\usepackage{amsmath, amscd, amsfonts, amssymb, amsthm, tikz, times}
\usepackage[T1]{fontenc}
\newtheorem{thm}{Theorem}[section]
\newtheorem{cor}[thm]{Corollary}
\newtheorem{lem}[thm]{Lemma}

\theoremstyle{definition}

\newtheorem*{remark}{Remark}

\newtheorem*{question}{Question}

\def\ed{\stackrel{d}{=}}
\def\bi{\begin{itemize}}
\def\ei{\end{itemize}}

\def\eps{\epsilon}

\def\be{\begin{align}}
\def\ee{\end{align}}
\def\bea{\begin{eqnarray}}
\def\eea{\end{eqnarray}}
\def\nn{\nonumber}
\def\ff{\infty}
\def\({\left(}
\def\){\right)}
\def\[{\left[}
\def\]{\right]}

\def\la{\left\langle}
\def\ra{\right\rangle}

\def\ER{Erd\H os-R\'enyi}

\def\t{\tilde}

\def\lam{\lambda}

\def\N{\mathbb{N}}
\def\R{\mathbb{R}}

\def\C{\mathbb{C}}

\def \P{\mathbf{P}}
\def \E{\mathbf{E}}
\def \O{\mathbf{O}}

\def \HH{{\cal H}}

\def\smath#1{\text{\scalebox{.8}{$#1$}}}

\def\sfrac#1#2{\smath{\frac{#1}{#2}}}

\def\ssqrt#1{\smath{\sqrt{#1}}}

\newcommand{\indic}{\mathbf{1}}

\newcommand{\Rom}[1]{\uppercase\expandafter{\romannumeral#1}}
\newcommand{\rom}[1]{\romannumeral#1\relax}

\begin{document}

\title{Delocalization and Limiting Spectral Distribution of Erd\H{o}s-R\'{e}nyi Graphs with Constant Expected Degree}
\author{{\normalsize{Paul Jung~,~~ Jaehun Lee~
}}\\[2ex]
\small{Department of Mathematical Sciences, Korea Advanced Institute of Science and Technology (KAIST),}\\[1ex]
\small{Daejeon 305-701, South Korea;\quad E-mail: pauljung@kaist.ac.kr, ljhiverson@kaist.ac.kr}\\[1ex]
}

\maketitle

\abstract{For fixed $\lambda>0$, it is known that Erd\H{o}s-R\'{e}nyi graphs $\{G(n,\sfrac{\lambda}{n}),n\in\N\}$, with edge-weights $\sfrac{1}{\sqrt{\lambda}}$, have a limiting spectral distribution, $\nu_{\lambda}$. As $\lambda\to\infty$,  $\{\nu_{\lambda}\}$ converges to the semicircle distribution. For large $\lambda$, we find an orthonormal eigenvector basis of $G(n,\sfrac{\lambda}{n})$ where most of the eigenvectors have small infinity norms as $n\to\infty$, providing a variant of an eigenvector delocalization result  of  Tran, Vu, and Wang (2013).}

\vspace{7mm}





\section{Introduction}
The spectral theory of graphs is important since many principal invariants of graphs are essentially related with their spectra. On the other hand, powerful tools used to investigate the spectrum of random matrices have been developed following the seminal work by Wigner \cite{Wigner1955}. In this paper, we study  a class of random matrices related to graphs, namely the adjacency matrices of Erd\H{o}s-R\'{e}nyi random graphs.

Let $G(n,p)$ be the Erd\H{o}s-R\'{e}nyi random graph with $n$ vertices and connection probability $p$. 
More precisely, letting $M_{n,p}$ denote the adjacency matrix of $G(n,p)$, for $i > j$ we independently set, 
\begin{align}
M_{n,p}(i,j)=\begin{cases}
	& 1 \;\; \text{with probability} \; p, \\
	& 0 \;\; \text{with probability} \; 1-p,
\end{cases}
\end{align}
and $M_{n,p}(i,j)=M_{n,p}(j,i)$ if $i < j$. Also, the graph has no loops, so $M_{n,p}(i,i)=0$ for all $i$. Note that $M_{n,p}$ is symmetric so its spectrum is real.

Recently, many outstanding results have been shown under the condition (with $p=p_n$)
\begin{align*}
\lim_{n\rightarrow\infty}np=\infty,
\end{align*}
in other words, under the condition that $G(n,p)$ has an expected degree, $np$, diverging with $n$. Under this condition, the spectral distribution of the scaled Erd\H{o}s-R\'{e}nyi ensemble
\begin{align*}
\frac{1}{\sqrt{np(1-p)}}M_{n,p}\;,\;\;n\in\N
\end{align*}
weakly converges to the standard semicircle distribution \cite{Tran2013}. Moreover, a local semicircle law holds \cite{Erdos2013}. Also, remarkably, all the $l^{2}$-normalized eigenvectors ``delocalize'' in term of their $l^{\infty}$-norm \cite{Erdos2013,Tran2013}. 

The situation is different if the expected degree is fixed.  If, for all $n$, we impose that $p=\lam/n$ for some fixed $\lam>0$, convergence to the semicircle law and delocalization do not hold \cite{Bauer2001, Bordenave2010, Zakharevich2005}. Let $\nu_{n,\lam}$ be the empirical spectral distribution of the scaled random adjacency matrix 
\begin{align}
\frac{1}{\sqrt{\lam}}M_{n,\lam/n}\begin{matrix}
\\.
\end{matrix}
\end{align}
As shown in \cite{Bauer2001, Bordenave2010, Zakharevich2005}, $\nu_{n,\lam}$ almost surely has a deterministic limiting distribution $\nu_{\lam}$ as $\lam$ goes to infinity; however, it is an open problem to find an explicit form for $\nu_{\lam}$, or even to give a characterization of its decomposition into pure-point, absolutely-continuous, and singular-continuous parts \cite{bordenave2013mean}. In \cite{Bauer2001}, Bauer and Golinelli analyzed $\nu_{\lam}$ using the moment method; we use the moment asymptotics given by their work as a starting point for this study. A numerical simulation is also given in \cite{Bauer2001}, and one can see that the numerical approximation of $\nu_{\lam}$ there, simulates the semicircle distribution as $\lam$ increases. 

\begin{thm}\label{main1}
	For each $\lam>0$, let $\nu_{n,\lam}$ be the empirical spectral distribution of $\sfrac{1}{\sqrt{\lam}}M_{n,\smath{\lam/n}}$. Let
	\begin{align*}
	\nu_{\lam}:=\lim_{n\rightarrow\infty}\nu_{n,\lam} \;\;\text{where the limit is in the weak sense.}
	\end{align*}
	Then, as $\lam$ goes to infinity, $\nu_{\lam}$ converges weakly to the standard semicircle distribution $\rho_{sc}$ where
	\begin{align*}
	\rho_{sc}(dx)=\sfrac{1}{2\pi}\sqrt{4-x^{2}}\indic_{\{|x|\le2\}}(dx).
	\end{align*}
\end{thm}
It was recently pointed out to us that the above result was proved in \cite{enriquez2016spectra}, nevertheless we provide two independent proofs of this fact since they are both different from the proof given in \cite{enriquez2016spectra}. These proofs are provided also for the sake of completeness, since the above result will play a crucial role in the proof of our main result, Theorem \ref{main2}.

Let us also remark that while the semicircle convergence results of \cite{Erdos2013, Tran2013} look similar to the above, there is a difference in the ``order of limits'': suppose $\{\lam_{m}\}$ is an expected degree sequence such that $\lim_{m\rightarrow\infty}\lam_{m}=\infty$. In \cite{Erdos2013, Tran2013}, a limiting ``diagonal'' spectral distribution sequence is considered,
\begin{align}\label{lim1}
\lim_{n\rightarrow\infty}\nu_{n,\lam_{n}},
\end{align}
whereas we are interested in the limit of limiting distributions $\{\nu_{\lam_{m}}\}$,
\begin{align}\label{lim2}
\lim_{m\rightarrow\infty}\nu_{\lam_{m}}=\lim_{m\rightarrow\infty}\lim_{n\rightarrow\infty}\nu_{n,\lam_{m}}.
\end{align}

 

In addition to results about the spectral distribution, another natural question is whether the $l^{2}$-normalized eigenvectors of $M_{n,\smath{\lam/n}}$ localize or delocalize. This question was raised, for example, by Dekel et al. \cite{Dekel2011}:

\begin{question}[Question 2 of \cite{Dekel2011}]\text{}\\
	\indent (\rom{1}) {\it Is it true that, almost surely, every unit eigenvector $u$ of $G(n,p)$ has $\lVert u \rVert_{\infty}=o(1)$?}\\
	\indent (\rom{2}) {\it Further, can we show that, almost surely, $\lVert u \rVert_{\infty}=n^{-\sfrac{1}{2}+o(1)}$?}
\end{question}

If the answer to (\rom{1}) is positive, we say that the unit eigenvectors {\bf delocalize}. Tao and Vu \cite{Tao2011} showed that (\rom{1}) and (\rom{2}) hold when $p=1/2$, which is of course independent of $n$. However, if $p=\lam/n$, it is easy to see that $G(n,p)$ almost surely has $\O(n)$ isolated vertices which persist in the limit. Thus, almost surely there exist at least $\O(n)$ eigenvectors such that their infinity norms are asymptotically 1, so delocalization fails. 

One can, however, obtain a weak form of delocalization as follows. For any $\eps>0$, one can choose $n$ and $\lam$ large enough so that most of the vectors in some $l^{2}$-normalized orthonormal basis have an infinity norm smaller than $\eps$. We need some notation in order to state this result more precisely. For any symmetric $n\times n$ matrix $H$, the eigenvalues of $H$ are denoted by $\{\Lambda_{i}(H)\}_{i=1}^{n}$. Without loss of generality, we suppose $$\Lambda_{1}(H)\le\Lambda_{2}(H)\le\cdots\le\Lambda_{n}(H)$$ throughout this paper. Since $H$ is symmetric, $H$ has an orthonormal basis $\{u_{i}(H)\}_{i=1}^{n}$ such that $u_{i}(H)$ is a unit eigenvector corresponding to $\Lambda_{i}(H)$.

\begin{thm}\label{main2}
	Let $\eps>0$. Using the above notation, define a subset $U(n,\lam,\eps)$ of $\{1,2,\cdots,n\}$ as follows,
	\begin{align}
	U(n,\lam,\eps):=\{i\in\smath{\{1,2,\cdots,n\}}:\lVert u_{i}(M_{n,\smath{\lam/n}}) \rVert_{\infty} < \eps \}.
	\end{align}
	Then, there exists an orthonormal basis $\{u_{i}(M_{n,\smath{\lam/n}})\}_{i=1}^{n}$ satisfying
	\begin{align*}
	\liminf_{\lam\rightarrow\infty}\liminf_{n\rightarrow\infty}\frac{|U(n,\lam,\eps)|}{n}= 1\;\;\;\text{almost surely}.
	\end{align*}
\end{thm}

	The strategy and main tools for proving the above are provided by  Theorem 1.16 in \cite{Tran2013} which we restate here for the reader's convenience.
	
\begin{thm}[Theorem 1.16 in \cite{Tran2013}]\label{TVWdeloc}
	Assume that the expected degree depends on $n$, i.e., $\lam=\lam_n$. Let $M_{n}:=M_{n,\lam_n/n}$. Suppose
	\begin{align}\label{prevcond}
	\lim_{n\rightarrow\infty}\frac{\lam_n}{\log{n}}=\infty.
	\end{align}
	Then there exists, a.s., an orthonormal eigenvector basis 
	\mbox{$
	\{u_{i}(M_{n}):i=1,2,\cdots,n\}
	$}
	such that
	\begin{align*}
	\lVert u_{i}(M_{n}) \rVert_{\infty} = o(1)
	\end{align*}
	for $1\le i \le n$.
\end{thm}

In fact, we also get a ``diagonalized convergence'' result as a corollary to Theorem \ref{main2}. The corollary should be viewed as a variant of the above Theorem \ref{TVWdeloc}. While the conclusion of the corollary is weaker than that of Theorem \ref{TVWdeloc}, the assumptions also allow for a broader class of sequences $\{\lambda_n\}$. This is one benefit of a priori considering the limiting behavior as two separate limits instead of one single diagonalized limit.

\begin{cor}\label{main3}
	Let $\lam=\lam_n$ depend on $n$ and set $M_{n}:=M_{n,\lam_n/n}$. Also, suppose	$\lim_{n\rightarrow\infty}\lam_n=\infty.$ Let $\eps>0$, and using the above notation, define $U'(n,\eps)$ by
	\begin{align}
	U'(n,\eps):=\{i\in\smath{\{1,2,\cdots,n\}}:\lVert u_{i}(M_{n}) \rVert_{\infty} < \eps \}.
	\end{align}
	Then, there exists a.s. an orthonormal eigenvector basis 
	$
	\{u_{i}(M_{n}):i=1,2,\cdots,n\}
	$	such that
	\begin{align*}
	\liminf_{n\rightarrow\infty}\frac{|U'(n,\eps)|}{n} = 1\;.
	\end{align*}
\end{cor}

The outline of the rest of this paper is as follows. In the next section (Section \ref{1stmain}), we give two proofs of Theorem \ref{main1} using respectively the moment method and the Stieltjes transform method. Section \ref{2ndmain} is devoted to the proofs of Theorem \ref{main2} and Corollary \ref{main3}.

\section{Convergence to the semicircle distribution}\label{1stmain}
As a preliminary to the two proofs, let us recall that the limiting distribution $\nu_{\lam}$ exists \cite{Bauer2001, Bordenave2010, Zakharevich2005}.  In particular, \cite{Bordenave2010} argues this via showing that the sequence of random graphs $\{G(n,\lam/n)\}_{n\in\N}$ converges, in the Benjamini-Schramm topology on rooted graphs, to a Galton-Watson tree with offspring distribution $\text{Pois}(\lam)$ (Poisson with intensity $\lam$). This fact will be useful to us in our second proof. Let us begin, however, with the classical moment method.

\subsection{Moment method proof}
\indent Fix $\lam>0$ and suppose $n\ge \lam$. Let $m_{ij}$
be the $(i,j)$ element of $M_{n,\smath{\lam/n}}$. A standard calculation in random matrix theory gives
\begin{align}\label{momen1}
\E\la \nu_{n,\lam},x^{k} \ra = \frac{1}{n\lam^{k/2}}\E \left[ \text{Tr}M_{n,\smath{\lam/n}}^{k} \right] = \frac{1}{n\lam^{k/2}} \sum_{1\le i_{1},\cdots,i_{k}\le n} \E \left[ m_{i_{1}i_{2}}m_{i_{2}i_{3}}\cdots m_{i_{k}i_{1}} \right].
\end{align}

We first obtain an asymptotic formula for $\E\la \nu_{n,\lam},x^{k} \ra$ using the method and terminology of \cite{Bauer2001}. 
If a $k$-tuple $(i_{1},i_{2},\cdots,i_{k})$ satisfies $i_{1}\neq i_{2}$, $i_{2}\neq i_{3}$, $\cdots$, $i_{k-1}\neq i_{k}$ and $i_{k}\neq i_{1}$, it is said to be {\bf admissible}. Non-admissible $k$-tuples do not contribute to the sum (\ref{momen1}) since $M_{n,\smath{\lam/n}}$ has vanishing diagonal entries. For each positive integer $j\le k$, define $W_{j}$ as the set of admissible $k$-tuple $(i_{1},i_{2},\cdots,i_{k})$ satisfying $|\{i_{1},i_{2},\cdots,i_{k}\}|=j$. The set $W$ of all admissible $k$-tuples is
\begin{align}
W:=\bigcup_{1\le j \le k}W_{j}.
\end{align}
\indent A $k$-tuple $(i_{1},i_{2},\cdots,i_{k})$ is called {\bf normalized} if it is admissible and $i_{j}>1$ implies that there exist $j'<j$ such that $i_{j'}=i_{j}-1$. Let $N_{j}$ be the set of normalized $k$-tuples $(i_{1},i_{2},\cdots,i_{k})$ such that $\{i_{1},i_{2},\cdots,i_{k}\}=\{1, 2, \cdots, j\}$. For $j\le n$, $\text{Per}(j,n)$ is defined to be the set of injective maps from $\{1,2,\cdots,j\}$ to $\{1,2,\cdots,n\}$. It is observed that, there is a one to one correspondence between $W_{j}$ and $\{(\omega,\sigma)| \omega\in N_{j}\;\text{and}\; \sigma\in\text{Per}(j,n)\}$. The set $N$ of all normalized $k$-tuples is expressed as
\begin{align}
N:=\bigcup_{1\le j \le k}N_{j}.
\end{align}
\indent In Eq. (\ref{momen1}), $m_{i_{1}i_{2}}m_{i_{2}i_{3}}\cdots m_{i_{k}i_{1}}$ can be identified with a closed walk along the graph given by the adjacency matrix $M_{n,\smath{\lam/n}}$. That is to say, $m_{i_{1}i_{2}}m_{i_{2}i_{3}}\cdots m_{i_{k}i_{1}}$ corresponds with the closed walk $i_{1}i_{2}\cdots i_{k}i_{1}$ (``closed'' means that it ends where it started). Let the sets of distinct edges and distinct vertices in the closed walk $i_{1}i_{2}\cdots i_{k}i_{1}$ corresponding to $k$-tuple $\omega=(i_{1},i_{2},\cdots,i_{k})$ be denoted by $E(\omega)$ and $V(\omega)$, respectively.  We denote an edge $e$ connecting the vertices with indices $i_{j}$ and $i_{j+1}$ by $e=i_{j}i_{j+1}$. Since $$m_{i_{1}i_{2}}m_{i_{2}i_{3}}\cdots m_{i_{k}i_{1}}=1$$ if and only if $m_{e}=1$ for all $e\in E(\omega)$,
\begin{align}
\frac{1}{n\lam^{k/2}} \sum_{1\le i_{1},\cdots,i_{k}\le n} \E \left[ m_{i_{1}i_{2}}m_{i_{2}i_{3}}\cdots m_{i_{k}i_{1}} \right] 
&=\nn \frac{1}{n\lam^{k/2}} \sum_{\omega\in W} \left(\sfrac{\lam}{n}\right)^{|E(\omega)|}\\
&=\frac{1}{n\lam^{k/2}} \sum_{\omega\in N} \left(\sfrac{\lam}{n}\right)^{|E(\omega)|} |\text{Per}(|V(\omega)|,n)| \label{tem1}
\end{align}

The moment method proof of Theorem \ref{main1} will follow from  Lemma \ref{momentlem} and Lemma \ref{momentmatchinglemma} below.
\begin{lem}\label{momentlem}
	For every positive integer $m$,
	\begin{align}\label{momenformula}
	\lim_{n\rightarrow\infty}\E\la \nu_{n,\lam},x^{k} \ra = \begin{cases}
	0 & k=2m-1 \\
	\frac{1}{m+1}\smath{\begin{pmatrix}2m \\ m \end{pmatrix}} + \O(\lam^{-1}) & k=2m
	\end{cases}
	\end{align}
	\begin{proof}
Let $\omega=(i_{1},i_{2},\cdots,i_{k})\in N$ and set $G(\omega)$ as the graph consisting of edges $E(\omega)$ and vertices $V(\omega)$.  We have $|E(\omega)|\ge |V(\omega)|-1$ since the graph $G(\omega)$ is connected.  
On the other hand, it is clear that in order to survive in the limit as $n\to\ff$ in (\ref{tem1}), one must have  $|V(\omega)|=|E(\omega)|+1$ because for any positive integer $j$
\begin{align*}
\lim_{n\rightarrow\infty}\frac{|\text{Per}(j,n)|}{n^{j}}=1 \;.
\end{align*}
In particular this implies that $G(\omega)$ must be a tree (rooted at $1$). 

Henceforth assume $|V(\omega)|=|E(\omega)|+1$. Then, $i_{1}i_{2}\cdots i_{k}i_{1}$ is a closed walk on a tree and so the multiplicity of every edge in the closed walk $i_{1}i_{2}\cdots i_{k}i_{1}$ is even. Thus, $2|E(\omega)|\le k$. 
Let $a_{l}$ be the number of normalized $k$-tuples $\omega$ such that $|E(\omega)|=l$ and $|V(\omega)|=|E(\omega)|+1$. 	In particular, if $k$ is odd, $a_{l}=0$ for all $1\le l \le k$ which proves the case $k=2m-1$ in \eqref{momenformula}. The $k=2m$ portion of (\ref{momenformula}) follows from (\ref{momen2}).
\begin{align}\label{momen2}
\lim_{n\rightarrow\infty} \frac{1}{n\lam^{k/2}} \sum_{\omega\in N} \left(\sfrac{\lam}{n}\right)^{|E(\omega)|} |\text{Per}(|V(\omega)|,n)|=\frac{1}{\lam^{ k/2 }} \sum_{l=1}^{\lfloor k/2 \rfloor} a_{l}\lam^{l}.
\end{align}
When $k=2m$, it is clear that $a_{m}$ is precisely the Catalan number, $C_m$, since the multiplicity of every edge in the closed walk $i_{1}i_{2}\cdots i_{k}i_{1}$ is exactly $2$.
	\end{proof}
\begin{remark}
	More precisely, when $k=2m$, one can easily check that
	\begin{align}
	\lim_{n\rightarrow\infty}\E\la \nu_{n,\lambda},x^{k} \ra = \frac{1}{m+1}\smath{\begin{pmatrix}2m \\ m \end{pmatrix}} + \sum_{l=1}^{m-1} a_{l}\lam^{l-m}.
	\end{align}
\end{remark}
\end{lem}

 \begin{lem}\label{momentmatchinglemma}
 	\begin{align}\label{momentmatching}
 	\lim_{n\rightarrow\infty}\E\la \nu_{n,\lam},x^{k} \ra = \la \nu_{\lam},x^{k} \ra
 	\end{align}
 	\begin{proof}
 	By Theorem 1 and Example 2 in \cite{Bordenave2010} (see also \cite[Thm 1.1]{Jung2017}), $\nu_{n,\lam}$ converges weakly to $\nu_\lam$ as $n\to\infty$. 
 		Thus, $\lim_{n\rightarrow\infty}\E\la \nu_{n,\lam},f \ra = \la \nu_{\lam},f \ra$ for any bounded continuous $f$, by dominated convergence. 
The lemma follows from a standard truncation argument. It is enough to consider the case for $k$ even because $\nu_{\lam}$ is symmetric (e.g., \cite[Thm 1.1]{Jung2017}). For $M>1$, define even functions $g_{M}$ with $g_{M}(x)=g_{M}(-x)$ by
 		\begin{align*}
 		&g_{M}(x)=\begin{cases}
 		1 & 0\le x\le M \\
 		0 & x \ge M+1 \\
 		-x+M+1 & M<x<M+1
 		\end{cases}
 		\end{align*}
so that
 		\begin{align*}
 		|\E\la \nu_{n,\lam},x^{2m} \ra - \E\la \nu_{n,\lam},x^{2m}g_{M} \ra|
 		\le \E\la \nu_{n,\lam},x^{2m}\indic_{|x|>M} \ra\le\frac{\E\la \nu_{n,\lam},x^{4m} \ra}{M^{2m}}
 		\end{align*}
Using the moment bound \eqref{momenformula}, take $n\to\infty$ then $M\to\infty$ to obtain (\ref{momentmatching}).
 	\end{proof}
%
 \end{lem}

 \hfill
 
 Recall that $\rho_{sc}$ is the standard semicircle distribution. It is easy to see that
 \begin{align}
 \lim_{\lam\rightarrow\infty}\la \nu_{\lam},x^{k} \ra=\la \rho_{sc},x^{k} \ra.
 \end{align}
Since $\rho_{sc}$ has bounded support, its moments characterize it uniquely, which implies that $\nu_{\lam}$ converges weakly to $\rho_{sc}$ (See Theorem 30.2 in \cite{Billingsley86}).\\
 \indent 

\subsection{Stieltjes transform proof}
\indent 
For later use, recall from \cite[pg. 225]{Reed1981} the notion of a  spectral measure $\nu_{\phi}$, of a self-adjoint operator $A$, associated to a unit vector $e_{\phi}$. Such a probability measure, $\nu_{\phi}$, can be defined by finding the unique measure satisfying
\begin{align}
\int_{\R}f(x)\nu_{\phi}(dx) =\la e_{\phi},f(A)e_{\phi} \ra
\end{align}
for all bounded, continuous $f$.

Using spectral theory and exchangeability, \cite{Bordenave2010} argued that the mean of the random measure $\nu_{n,\lam}$ can be regarded as the expected spectral measure at vertex $1$ (or any other fixed vertex) of the Erd\H{o}s-R\'{e}nyi graph $G(n,\lam/n)$ (with weights $1/\sqrt{\lam}$ on the edges). Moreover, the limiting deterministic measure $\nu_\lam$ is the expected spectral measure associated to the root of a Galton-Watson tree with offspring distribution  $\text{Pois}(\lam)$ and weights $1/\sqrt{\lam}$, which is the limit of $\{G(n,\lam/n)\}$ with weighted edges in the Benjamini-Schramm topology (see also \cite{Bordenave2011,Jung2017}). The adjacency operator $\sfrac{1}{\sqrt{\lam}}M_{\ff}^{(\lam)}$ of the limiting graph is self-adjoint (\cite[Lemma 5.2]{Jung2017}) and its resolvent $R^{(\lam)}$ is well-defined. Letting $\phi$ denote the root of the tree and $e_\phi$ denote the root vector, i.e. a Kronecker-delta function at the root,
define the random variable $$R^{(\lam)}_{\phi,\phi}(z):=\la e_{\phi},\left(\sfrac{1}{\sqrt{\lam}}M_{\ff}^{(\lam)}-zI\right)^{-1}e_{\phi} \ra$$
where the domain of $z$ is $\C\backslash\R$.

Let $S_{\lam}$ be the Stieltjes transform of the limiting distribution $\nu_{\lam}$. According to \cite[Thm 2]{Bordenave2010},
\begin{align}\label{rde}
R^{(\lam)}_{\phi,\phi}(z)\ed -\left[ \frac{1}{z+\tfrac{1}{\lam}\sum_{k=1}^{\text{Pois}(\lam)}R^{(\lam)}_{k,k}(z)}\right]
\end{align}
where $(R^{(\lam)}_{k,k}(z))_{k\in\N}$ is an i.i.d. sequence with the same distribution as $R^{(\lam)}_{\phi,\phi}(z)$ and $\text{Pois}(\lam)$ is a Poisson random variable independent from $(R^{(\lam)}_{k,k}(z))_{k\in\N}$. Thus, $$S_{\lam}(z)=\E R^{(\lam)}_{\phi,\phi}(z).$$

The strategy of the proof is to show that 
$$S(z):=\lim_{\lam\rightarrow\infty}S_{\lam}(z)$$ exists for all $z\in\C\backslash\R$ and
 satisfies the self-consistent equation,
\begin{align}\label{key1}
S(z)=-\frac{1}{z+S(z)}
\end{align}
implying that $S(z)=-\frac12(z-\sqrt{z^2-4})$ by choosing the solution of \eqref{key1} such that the imaginary parts of $S(z)$ and $z$ are the same.
By the Stieltjes inversion formula, $\nu_{\lam}$ converges weakly to $\rho_{sc}$, the standard semicircle law, as $\lam\to\ff$. 

Let us now carry out the above strategy. Define $Y_\lam$ and $f_{\lam}$ as follows.
$$Y_\lam:=\frac{1}{\lam}\sum_{k=1}^{\text{Pois}(\lam)}R^{(\lam)}_{k,k}(z) \;\; \text{and} \;\; f_{\lam}(\theta):=\E \exp\big(i\theta R^{(\lam)}_{\phi,\phi}(z)\big)$$
so that
\begin{align}\label{calculation}
\E\exp(i\theta Y_\lam)&=\E\big[ \big\{f_{\lam}(\tfrac{\theta}{\lam})\big\}^{\text{Pois}(\lam)}\big]=\E\big[ \exp[\text{Pois}(\lam)\log f_{\lam}(\tfrac{\theta}{\lam})] \big] \nn\\
&=\exp\big[\lam\big( e^{\log f_{\lam}(\sfrac{\theta}{\lam})}-1 \big)\big]=\exp\big[\lam\big(f_{\lam}(\smath{\sfrac{\theta}{\lam}})-1 \big)\big]\nn\\
&=\exp\big[\lam\big(\tfrac{i\theta}{\lam}\E R^{(\lam)}_{\phi,\phi}(z)+o(\tfrac{1}{\lam})\big)\big]
\end{align}
where $o(\tfrac{1}{\lam})$ depends on $\theta$. The last equality in (\ref{calculation}) comes from the Taylor expansion of the characteristic function $f$ which is possible since we have the a.s. bound
\begin{align}\label{trivialbound}
 \left|R^{(\lam)}_{\phi,\phi}(z)\right| \le \int \left|\frac{1}{x-z}\right| d\nu_{\phi}\le\frac{1}{|\text{Im}(z)|}.
\end{align}

Choose a subsequence $\lam_n\to\infty$ such that a limit $S(z)$ exists. Eq. \eqref{calculation} tells us that
\begin{align}\label{eq1}
\frac{1}{\lam_n}\sum_{k=1}^{\text{Pois}(\lam_n)}R^{(\lam_n)}_{k,k}(z)\overset{pr}{\longrightarrow} S(z) \;\; \text{as}\; n\rightarrow\infty,
\end{align}
by convergence of the characteristic functions of $\{Y_{\lam_n}, \lam_n>0\}$, and the fact that the limit is a constant.  
Next, suppose without loss of generality that $z\in\C^+$. Then
$$\Im(S(z)+z)\ge \Im(z)>0$$ which implies $S(z)\neq z$. 
 By the continuous mapping theorem,
\begin{align}\label{den}
-\frac{1}{z+\frac{1}{\lam_n}\sum_{k=1}^{\text{Pois}(\lam_n)}R^{(\lam_n)}_{k,k}(z)}\overset{pr}{\longrightarrow}-\frac{1}{z+S(z)}.
\end{align}
 By \eqref{rde}, the left-hand side above has the same distribution as $R^{(\lam_n)}_{\phi,\phi}(z)$ which by \eqref{trivialbound} is bounded for any fixed  $z\in \C\backslash\R$. Thus
\begin{align}\label{gap1}
\lim_{n\rightarrow\infty} S_{\lam_n}(z) = \lim_{n\rightarrow\infty}\E R^{(\lam_n)}_{\phi,\phi}(z) = -\frac{1}{z+S(z)}.
\end{align} 
Therefore $S(z)$ satisfies (\ref{key1}) and must be the Stieltjes transform of the semicircle law. The proof follows since the measures $\{\nu_{\lam}\}$ are tight 
, while the above argument shows that there is a unique limit point.

\section{Delocalization}\label{2ndmain}
Recall that $\{\Lambda_{i}(H)\}_{i=1}^{n}$ and $\{u_{i}(H)\}_{i=1}^{n}$ denote the eigenvalues and eigenvectors of a symmetric $n\times n$ matrix $H$, respectively. 
We begin with several lemmas, the first of which is Eq. (5.8) in \cite{erdHos2009semicircle}.  We state the version from \cite[Lemma 41]{Tao2011}:

\begin{lem}[Lemma 41 in \cite{Tao2011})]\text{}\label{schur}\\
	Let
	\begin{align}
	H=\begin{pmatrix}
	a & X^{T} \\
	X & \t{H}
	\end{pmatrix}
	\end{align}
	be an $n\times n$ symmetric matrix for some $a\in\R$ and $X\in \R^{n-1}$, and let $\smath{\begin{pmatrix}x \\ v \end{pmatrix}}$ be the unit eigenvector with eigenvalue $\Lambda_{i}(H)$ where $x\in\R$ and $v\in\R^{n-1}$. Assume that none of the eigenvalues of $\t{H}$ are equal to $\Lambda_{i}(H)$. Then,
	\begin{align}
	|x|^{2}=\frac{1}{1+\sum_{j=1}^{n-1}(\Lambda_{j}(\t{H})-\Lambda_{i}(H))^{-2}|\langle u_{j}(\t{H}),X\rangle|^{2}}
	\end{align}
	where $\la\cdot,\cdot\ra$ denotes the inner product between vectors.
\end{lem}

The second lemma is a consequence of Talagrand's inequality that was proved in Lemma 68 of \cite{Tao2011}. We state the version from \cite[Lemma 3.4]{Tran2013}:

\begin{lem}[Lemma 3.4 in \cite{Tran2013}]\label{concentdegree}
	Let $Y=(y_{1}, \cdots, y_{n})\in\R^{n}$ be a random vector whose coordinates are i.i.d. centered random variables which are a.s. {bounded in absolute value by 1} and have variance $\sigma^{2}$. Let $\HH$ be a subspace of dimension $k$ and $\pi_{\HH}$ the orthogonal projection onto $\HH$. Then,
	\begin{align}
	\P \left( \left|\;\lVert \pi_{\HH}(Y) \rVert - \sigma\sqrt{k}\;\right|\ge t \right) \le 10\exp(-t^{2}/4)
	\end{align}
	where $\lVert \cdot \rVert$ is the Euclidean norm.
\end{lem}

Let $N_{n}$ be a symmetric $n\times n$ matrix whose upper triangular elements are independent standard normal variables $N(i,j)$.
Note that even though the perturbed matrix elements are unbounded, we have that
\begin{align*}
\P[|N(i,j)|>\sqrt{n}]\le Ce^{-\frac{n}{2}}.
\end{align*}
As $\sum_{n}n^2e^{-\frac{n}{2}} < \ff $, by Borel-Cantelli we have that $|N(i,j)|\le \sqrt{n}$ a.s. for all $1\le i,j,\le n$ and all $n$ large enough. This will allow us to use Lemma \ref{concentdegree} later on.

Assume that $N_{n}$ is also independent from $M_{n,\smath{\lam/n}}$. Let $\{\delta(n)\}_{n\in\N}$ be a sequence of positive numbers satisfying
\begin{align}\label{decay}
\delta(n)=o(n^{-1/2}).
\end{align}
Denote the scaled adjacency matrix and a perturbed version of it as follows:
\begin{align}
A_{n,\lam}&:=\frac{1}{\sqrt{\lam}}M_{n,\smath{\lam/n}} \label{start}\;,\\
B_{n,\lam}&:=A_{n,\lam}+\delta(n)N_{n}\;.
\end{align}
The reason for introducing the perturbed matrix is that it almost surely has a simple spectrum (see \cite[Exercise 1.3.10]{tao2012topics}):
\begin{align}\label{simplespectrum}
\Lambda_{1}(B_{n,\lam})<\Lambda_{2}(B_{n,\lam})<\cdots<\Lambda_{n}(B_{n,\lam})\;\;\text{almost surely.}
\end{align}

Write $B_{n,\lam}$ in the following matrix form:
\begin{align}\label{def:B}
B_{n,\lam}&=\begin{pmatrix}a & X^{T} \\ X & \t{B}_{n,\lam} \end{pmatrix} \;\; \text{where}\; a\in\R \;\text{and}\; X\in\R^{n-1}.
\end{align}
Then,
\begin{align}\label{goodcondition2}
\{\Lambda_{i}(B_{n,\lam}):i=1,2,\cdots,n\} \cap \{\Lambda_{i}(\t{B}_{n,\lam}):i=1,2,\cdots,n-1\} = \emptyset  \;\;\text{almost surely}
\end{align}
by \eqref{simplespectrum} and the Cauchy interlacing principle. Note that \eqref{goodcondition2} allows us to use Lemma \ref{schur}. Our third preliminary lemma bounds the effect of the above perturbation on infinity norms of eigenvectors:
\begin{lem}[Lemma 3.1 in \cite{Tran2013}]\label{goodtool} \text{}\newline
	Recall that $B_{n,\lam}$ is defined as the perturbation $A_{n,\lam}+\delta(n)N_{n}$. There exists an orthonormal basis of eigenvectors $\{u_{i}(A_{n,\lam})\}_{i=1}^{n}$ such that, for every $1\le i \le n$,
	\begin{align}\label{goodcondition}
	\lVert u_{i}(A_{n,\lam}) \rVert_{\infty} \le \lVert u_{i}(B_{n,\lam}) \rVert_{\infty} + \alpha(n)
	\end{align}
	where $\alpha(n)\to 0$ as $n\to\infty$, and $\alpha(n)$ can be chosen to be arbitrarily small depending only on $\delta(n)$.
\end{lem}

Henceforth assume $u_{i}(A_{n,\lam})=u_{i}(M_{n,\smath{\lam/n}})$ for all $i$ and $n$ and that the orthonormal basis $\{u_{i}(A_{n,\lam})\}_{i=1}^{n}$ satisfies (\ref{goodcondition}).


\begin{lem}\label{convlem}	Let $\t{\mu}_{n,\lam}$ be the empirical spectral distribution of $\t{B}_{n,\lam}$. Then, for $a<b$, and $\delta(n)$ satisfying \eqref{decay},
	\begin{align}\label{1sttarget}
	\limsup_{\lam\rightarrow\infty}\limsup_{n\rightarrow\infty} \left|\nu_{\lam}([a,b]) - \t{\mu}_{n,\lam}([a,b])   \right| = 0, \;\text{almost surely}.
	\end{align}
\end{lem}
The above lemma follows simply from Theorem \ref{main1} and Weyl's inequality; however, for completeness, we provide an explicit proof in Appendix \ref{appen1}.

\subsection{Proof of Theorem \ref{main2}}

By Lemma \ref{goodtool}, it is sufficient to show that the conclusion of Theorem \ref{main2} holds when $u_{i}(M_{n,\smath{\lam/n}})$ is replaced with $u_{i}(B_{n,\lam})$. Let $\tilde{U}$ be defined by
\begin{align}
\tilde{U}(n,\lam,\eps):=\{i\in\smath{\{1,2,\cdots,n\}}:\lVert u_{i}(B_{n,\lam}) \rVert_{\infty} < \eps/2 \}.
\end{align}
Our goal is to prove
\begin{align*}
\liminf_{\lam\rightarrow\infty}\liminf_{n\rightarrow\infty}\frac{|\tilde{U}(n,\lam,\eps)|}{n}=1\;\;\;\text{almost surely}.
\end{align*}
For this, it suffices to  show
\begin{align}\label{eq: first part}
\limsup_{\lam\rightarrow\infty}\limsup_{n\rightarrow\infty}\lVert u_{i}(B_{n,\lam}) \rVert_{\infty} < \frac{\eps}{2} \;\;\;\text{almost surely,}
\end{align}
uniformly for all $\Lambda_{i}(B_{n,\lam})\in[-2,2]$,
since by Lemma \ref{convlem}
\begin{align}\label{eq: second part}
\liminf_{\lam\rightarrow\infty}\liminf_{n\rightarrow\infty} \frac{|\{i\in\smath{\{1,2,\cdots,n\}}:\Lambda_{i}(B_{n,\lam})\in[-2,2]\}|}{n}=1 \;\;\;\text{almost surely}.
\end{align}

By (\ref{goodcondition2}), we can apply Lemma \ref{schur} to get
\begin{align}\label{eq:x bound}
|x|^{2}=\frac{1}{1+\sum_{j=1}^{n-1}(\Lambda_{j}(\t{B}_{n,\lam})-\Lambda_{i}(B_{n,\lam}))^{-2}\left|\la u_{j}(\t{B}_{n,\lam}),X\ra\right|^{2}}
\end{align}
where $x=u_i(1)$ is the first coordinate of $u_{i}(B_{n,\lam})$. A similar bound holds for any other coordinate $u_i(k)$ of $u_{i}(B_{n,\lam})$ by replacing $\t{B}_{n,\lam}$ with an appropriate submatrix. Thus, we will see that it suffices to find an upper bound of $|x|^{2}$, with high enough probability, in order to get an upper bound for $\lVert u_{i}(B_{n,\lam}) \rVert_{\infty}$, uniformly in $i$ with high probability. 

Let $Q$ be a positive integer and set 
\begin{align}\label{def:l}
l:=4/Q.
\end{align}
Choose  $Q$ large enough so that $Q\ge 5$ and
\begin{align}\label{eq: eps over 2}
\frac{1}{1+1/(\pi\sqrt{3l})}<\frac{\eps^{2}}{4}.
\end{align}
We fix this value of Q (thus fixing $l$) henceforth and note that they only depend on $\eps$. 

Partition the interval $[-2,2]$ into $\{[a_{q},a_{q+1}]\}_{q=1}^{Q}$ so that $a_{1}=-2$, $a_{Q+1}=2$ and $a_{q+1}-a_{q}=l$ for every $q$. Suppose now that $\Lambda_{i}(B_{n,\lam})\in[-2,2]$ so that there is a $q_i$ such that $\Lambda_{i}(B_{n,\lam})\in[a_{q_i},a_{q_i+1}]$.
Define a subset $J$ of $\{1,2,\cdots, n-1\}$ as 
\begin{align}
J(n,q):=\{j\in\smath{\{1,2,\cdots, n-1\}:\Lambda_{j}(\t{B}_{n,\lam}) \in [a_{q}, a_{q+1}]}\}
\end{align}
and define $\HH=\HH(n,q):=\text{span}_{j\in J(n,q)}\{u_{j}(\t{B}_{n,\lam})\}$. Let $\pi_{\HH}$ be the orthogonal projection onto $\HH$. For $q_i$ as above
\begin{align}
J(n,q_i)\subset \{j\in\smath{\{1,2,\cdots, n-1\}:|\Lambda_{j}(\t{B}_{n,\lam})-\Lambda_{i}(B_{n,\lam})|\le l}\}.
\end{align}
Therefore {for $\HH_i=\HH(n,q_i)$}, we get the inequality
\begin{align}\label{eq:lower bound}
\nn\sum_{j=1}^{n-1}(\Lambda_{j}(\t{B}_{n,\lam})-\Lambda_{i}(B_{n,\lam}))^{-2}\left|\la u_{j}(\t{B}_{n,\lam}),X\ra\right|^{2}&\ge \frac{1}{l^{2}}\sum_{j\in J(n,q_i)}\left|\la u_{j}(\t{B}_{n,\lam}),X\ra\right|^{2}\\
&=\frac{1}{l^{2}}\lVert \pi_{\HH_i}(X) \rVert^{2}.
\end{align}

Now, define a random vector $Y$ from the vector $X$ which is as in \eqref{def:B}
{\begin{align}
Y:=X-\frac{\ssqrt{\lam}}{n}\indic(n-1) \;\; \text{where}\;\;\indic(n)=(1,1,\cdots,1)\in\R^{n}.
\end{align}
Let $\HH'$ be the orthogonal complement of $\text{span}\{\indic(n-1)\}$.} Then, for generic $q$ and  $\HH=\HH(n,q)$,
\begin{align}
\lVert  \pi_{\HH}(X) \rVert \ge \lVert  \pi_{\HH\cap \HH'}(X) \rVert =  \lVert  \pi_{\HH\cap \HH'}(Y) \rVert.
\end{align}
Observe in particular, that $
\text{dim}(\HH\cap \HH')\ge\text{dim}(\HH)-1.$
Since $\t{B}_{n,\lam}$ is independent of $Y$, Lemma \ref{concentdegree} can be applied with
\begin{align}
t=t(n)=\sqrt{\sqrt{n}\cdot\log{n}}
\end{align}
after conditioning on $\t{B}_{n,\lam}$, and also after normalizing $Y$ so that $\sigma=1$. Thus with probability at least $1-10\exp(-(\sqrt{n}\cdot\log{n})/4)$,
\begin{align}\label{concentration4}
\left\lVert \pi_{\HH\cap \HH'}\left(\left(\sfrac{1-\lam/n}{n}+\delta^{2}\right)^{-\frac{1}{2}}\cdot Y\right) \right\rVert  \ge \sqrt{|J(n,q)|-1}-\ssqrt{\sqrt{n}\cdot\log{n}}.
\end{align}
The Borel-Cantelli lemma implies that the inequality (\ref{concentration4}) holds almost surely for  large $n$ and for every subspace $\HH(n,q)$ with $q\in[1,Q]$, {so in particular it holds for $\HH(n,q_i)$.} Plugging \eqref{concentration4} into \eqref{eq:lower bound}, and
recalling that $\delta(n)=o(n^{-1/2})$, we have, almost surely,
{\begin{align}\label{eq: 52}
\liminf_{n\rightarrow\infty}\sum_{j=1}^{n-1}(\Lambda_{j}(\t{B}_{n,\lam})-\Lambda_{i}(B_{n,\lam}))^{-2}\left|\la u_{j}(\t{B}_{n,\lam}),X\ra\right|^{2} \ge \liminf_{n\rightarrow\infty}\min_{q\in\{1,2,\cdots,Q\}}\frac{|J(n,q_i)|}{n\cdot l^{2}}
\end{align}}

Recall that $\t{\mu}_{n,\lam}$ is the empirical spectral distribution of $\t{B}_{n,\lam}$. Fix $q\in\{1,2,\cdots,Q\}$ and note that
\begin{align}
\frac{|J(n,q)|}{n-1}=\t{\mu}_{n,\lam}([a_{q},a_{q+1}]). \label{end}
\end{align}
Applying Theorem \ref{main1} and Lemma \ref{convlem} to
\begin{align}\label{keyineq2}
\t{\mu}_{n,\lam}([a_{q},a_{q+1}]) \ge&\;\; \rho_{sc}([a_{q},a_{q+1}])- \left| \rho_{sc}([a_{q},a_{q+1}]) - \nu_{\lam}([a_{q},a_{q+1}]) \right|\nn\\
&-\left|\nu_{\lam}([a_{q},a_{q+1}]) - \t{\mu}_{n,\lam}([a_{q},a_{q+1}])   \right|
\end{align}
and using the calculation in Appendix \ref{appen2}, we have almost surely,
\begin{align}
\liminf_{\lam\rightarrow\infty}\liminf_{n\rightarrow\infty}\min_{q\in\{1,2,\cdots,Q\}}\frac{|J(n,q)|}{n-1} &\ge \min_{q\in\{1,2,\cdots,Q\}}\rho_{sc}\big([a_{q},a_{q+1}]\big)\nn\\
&\ge\frac{l^{3/2}}{\pi\sqrt{3}}\;.\label{tem2}
\end{align}
Combining \eqref{eq:x bound}, \eqref{eq: eps over 2}, \eqref{eq: 52} and \eqref{tem2} we get that $|x|<\eps/2$ almost surely for large $n$ and large $\lam$ under the assumption $\Lambda_{i}(B_{n,\lam})\in[-2,2]$. 

Finally, recall that a relation similar to \eqref{eq:x bound} holds for any other coordinate $u_i(k)$ of $u_i(\t{B}_{n,\lam})$ and so using a union bound over $k\in\{1,\ldots,n\}$, and noting that $\sum_n 10n\exp(-(\sqrt{n}\cdot\log{n})/4)<\infty$ (in order to invoke the Borel-Cantelli lemma for a union of probabilities), we obtain \eqref{eq: first part}. This completes the proof.


\subsection{Proof of Corollary \ref{main3}}

From now on, we let the expected degree depend on $n$, i.e., $\lam=\lam_n$. 
Recall that in contrast to Theorem \ref{TVWdeloc} where the growth condition (\ref{prevcond}) is required, we consider the more general case where
\begin{align}
\lim_{n\rightarrow\infty}\lam_{n}=\infty.
\end{align}
 Recall that $M_{n}:=M_{n,\lam_n/n}$. Also, let $\nu_{n}:=\nu_{n,\lam_n}$.
According to Theorem 1.3 in \cite{Tran2013}, the empirical spectral measure $\nu_{n}$ weakly converges to the standard semicircle distribution $\rho_{sc}$ as $n$ goes to infinity. We can use the same argument as in the proof of Theorem \ref{main2} up until (\ref{end}). After that, set $\mu_{n}:=\mu_{n,\lam_{n}}$ and $\tilde{\mu}_{n}:=\tilde{\mu}_{n,\lam_{n}}$ and use the following inequality instead of (\ref{keyineq2}):
\begin{align}
\t{\mu}_{n}([a_{q},a_{q+1}]) \ge&\;\; \rho_{sc}([a_{q},a_{q+1}])- \left| \rho_{sc}([a_{q},a_{q+1}]) - \t{\mu}_{n}([a_{q},a_{q+1}]) \right|.
\end{align}
By the absolute continuity of $\rho_{sc}$, and the argument in Appendix \ref{appen1}, we have
\begin{align}
\limsup_{n\rightarrow\infty}\left| \rho_{sc}([a_{q},a_{q+1}]) - \t{\mu}_{n}([a_{q},a_{q+1}]) \right| = 0.
\end{align}
Consequently,
\begin{align}
\liminf_{n\rightarrow\infty}\min_{q\in\{1,2,\cdots,Q\}}\frac{|J(n,q)|}{n\cdot l^{2}} \ge \min_{q\in\{1,2,\cdots,Q\}}\frac{\rho_{sc}\big([a_{q},a_{q+1}]\big)}{l^{2}}\ge \frac{1}{\pi\sqrt{3l}}\;\begin{matrix}
\\.
\end{matrix}
\end{align}
Since
$\liminf_{n\rightarrow\infty}\mu_{n}([-2,2])=1$, the result follows.

While Corollary \ref{main3} has the advantage of holding without any growth rate condition on $\lambda_n$, it has the drawback that it give no information about the infinity norms of eigenvectors corresponding to the eigenvalues outside of $[-2,2]$. Note that $[-2,2]$ corresponds to the support of the standard semicircle law.




\hfill

\appendix
\section{Some additional tools}
\subsection{Proof of Lemma \ref{convlem}}\label{appen1}
{With some abuse of notation, write $\t{B}_{n,\lam}=A_{n-1,\lam}+\delta(n)N_{n-1}$. Then, Weyl's theorem implies
\begin{align*}\label{sameglobal}
|\Lambda_{i}(\t{B}_{n,\lam})-\Lambda_{i}(A_{n-1,\lam})|\le\delta(n)\lVert N_{n-1} \rVert_{\text{op}}=O(\delta(n)\sqrt{n})\;\;\text{for all}\;1\le i \le n
\end{align*}
since $\{(n)^{-1/2}N_{n}\}_{n\in\N}$ is a Wigner ensemble with moments of all order. Using \eqref{decay}, Weyl's inequality, and the Cauchy interlacing theorem, there is a sequence $\lim_{n\rightarrow\infty}\zeta_{n}=0$ such that
\begin{align}
\nu_{n-1,\lam}([a+\zeta_{n},b-\zeta_{n}]) \le \t{\mu}_{n,\lam}([a,b]) \le \nu_{n-1,\lam}([a-\zeta_{n},b+\zeta_{n}]).
\end{align}}
Note that
\begin{align*}
\limsup_{n\rightarrow\infty}\nu_{n,\lam}([a-\zeta_{n},b+\zeta_{n}]) \le \lim_{\xi\downarrow0}\limsup_{n\rightarrow\infty}\nu_{n,\lam}([a-\xi,b+\xi]).
\end{align*}
For fixed $\xi>0$, choose continuous functions $f_{\xi}$ and $g_{\xi}$ which converge pointwise to $\indic_{[a,b]}$, as $\xi\to\ff$, and which satisfy
\begin{align*}
0\le f_{\xi}\le \indic_{[a-\xi,b+\xi]} \le g_{\xi} \le 1.
\end{align*}
Then, almost surely,
\begin{align*}
\int f_{\xi}\;d\nu_{\lam}=\lim_{n\rightarrow\infty}\int f_{\xi}\;d\nu_{n,\lam}\le \limsup_{n\rightarrow\infty}\nu_{n,\lam}([a-\xi,b+\xi])\le\lim_{n\rightarrow\infty}\int g_{\xi}\;d\nu_{n,\lam}=\int g_{\xi}\;d\nu_{\lam}\;.
\end{align*}
We deduce that
\begin{align}
\limsup_{n\rightarrow\infty}\left|\nu_{\lam}([a,b]) - \t{\mu}_{n,\lam}([a,b])   \right| = \nu_{\lam}(\{a\})+\nu_{\lam}(\{b\})\;.
\end{align}
Finally by Theorem \ref{main1}, both $\nu_{\lam}(\{a\})$ and $\nu_{\lam}(\{b\})$ go to 0 as $\lam\to\infty$. 
\subsection{A simple bound for the semicircle edge}\label{appen2}
Recall from \eqref{def:l} that $l=4/Q$. Here we will show that when $[-2,2]$ is partitioned into $Q$ equal parts,
\begin{align}
\min_{q\in\{1,2,\cdots,Q\}}\rho_{sc}\big([a_{q},a_{q+1}]\big)\ge\frac{l^{3/2}}{\pi\sqrt{3}}\begin{matrix}
\\.
\end{matrix}
\end{align}
Observe that
\begin{align*}
\min_{q\in\{1,2,\cdots,Q\}}\rho_{sc}\big([a_{q},a_{q+1}]\big)=\frac{1}{2\pi}\int_{2-l}^{2}\sqrt{4-x^{2}}\;dx.
\end{align*}
Since $4-x^{2}=(2+x)(2-x)$, we have for $l<1$,
\begin{align*}
\frac{1}{2\pi}\int_{2-l}^{2}\sqrt{4-x^{2}}\;dx \ge \frac{\sqrt{3}}{2\pi}\int_{0}^{l}\sqrt{x}\;dx=\frac{l^{3/2}}{\pi\sqrt{3}}\begin{matrix}
\\.
\end{matrix}
\end{align*}

\paragraph{Acknowledgments} Research supported by NRF grant 2017001952 (Republic of Korea). We thank Terence Tao for pointing out reference \cite{Tran2013} to us, as well as an anonymous referee for helpful comments.

\vspace{5mm}

\end{document}